\documentclass[12pt]{article}
\usepackage{amsmath,amsthm,amssymb,latexsym,color}
\usepackage{bbm}

\usepackage{graphicx}
\usepackage{tikz}
\usetikzlibrary{decorations.pathreplacing,shapes.misc}

\usepackage{pgfplots}

\usepackage{hyperref}

\date{}

\theoremstyle{plain}

\newtheorem*{theorem*}{Theorem}

\newtheorem{theor}{Theorem}[section]

\newtheorem*{theor*}{Theorem}
\newtheorem*{conj*}{Conjecture}

\newtheorem{thm}{Theorem}[section]

\newtheorem{prop}[theor]{Proposition}
\newtheorem{defin}[theor]{Definition}

\newtheorem{lemma}[theor]{Lemma}
\newtheorem{rem}[theor]{Remark}
\theoremstyle{remark}

\def\Bo{{\mathcal B}}
\def\Be{\Bo _{\eps, d}}
\def\Bep{\wt \Bo _{\eps, d}}
\def\Bem{\Bo _{\eps, m}}
\def\Bepm{\wt \Bo _{\eps, m}}
\def\Rc{{\mathcal R}}
\def\Nc{{\mathcal N}}
\def\Mc{{\mathcal M}}
\def\Ec{{\mathcal E}}
\def\R{{\mathbb R}}
\def\N{{\mathbb N}}

\def\Prob{{\mathbb P}}

\def\r{\right}
\def\lf{\lfloor}
\def\rf{\rfloor}

\def\wt{\widetilde}

\def\d{{\rm disp}}

\newcommand{\eps}{\varepsilon}

\def\qand{\quad \mbox{ and } \quad}

\title{A remark on the minimal dispersion}
\author{A. E. Litvak}

\newcommand\address{\noindent\leavevmode

\medskip
\noindent
Alexander  Litvak\\
Dept.~of Math.~and Stat.~Sciences,\\
University of Alberta, \\
Edmonton, AB, Canada, T6G 2G1.\\
\texttt{\small
e-mail:  alitvak@ualberta.ca}
}

\begin{document}

\maketitle

\begin{abstract}
We improve known upper bounds for the minimal dispersion of a point set
in the unit cube and its inverse in both the periodic and non-periodic settings.
Some of our bounds are sharp up to logarithmic factors.
\end{abstract}

\bigskip

{\small
\noindent{\bf AMS 2010 Classification:}
primary: 52B55, 52A23;
secondary: 68Q25, 65Y20.\\
\noindent
{\bf Keywords:} complexity, dispersion, largest empty box, torus






\section{Introduction and main results}
\label{intro}

In this note we deal with the minimal dispersion of a point set in the unit cube.
The dispersion of a point set in the $d$-dimensional unit cube $[0,1]^d$ is defined as the maximal volume
of an axis parallel box in the cube which does not contain any point from the set.
Then the minimal dispersion is a function of two variables, $n$ and $d$, which minimizes
the dispersion over all possible choices of $n$ points. Such definition was introduced
in \cite{RT}  modifying a  notion from \cite{Hl}. Due to important applications and due
to the fact that the problem is very interesting by itself, it has attracted a considerable attention
in recent years. We refer to \cite{AHR, DJ, HKKR, Rud, Sos,  MU, UV, UV1} and references therein for the
history of the problem and its relation to other branches as well as for the best known bounds
(see also \cite{Kr, VT, MU-lat} for the dispersion of certain sets).
We improve known upper bounds for the minimal dispersion and its inverse function.
We will also consider the minimal dispersion on the torus  and discuss the sharpness of our results.
We would like to emphasize that we look at the dispersion as at a function of two variables,
without trying to fix one of the variables. Instead, we consider both variables growing to infinity and
our bounds depend on the relations between variables. The main novelty in our proof is a better
construction of a family of axis parallel boxes  (periodic or non-periodic)
needed to be checked for a random choice of points.
It seems that our construction also leads to better bounds for  recently introduced in \cite{HPUV} $k$-dispersion
(where, given set of $n$ points, one allows axis parallel boxes  to have inside  at most $k$ points from this set),
 but we do not pursue this direction.

\subsection{Notations}

We start with notations.
Given a measurable set $A\subset \R^d$, we denote its $d$-dimensional volume by $|A|$.
We also use the same notation $|M|$ for the cardinality of a finite set $M$
(it always will be clear from the context what $|\cdot|$ means).
By $\Rc _d$
we denote the set of all axis parallel boxes contained in the cube
$Q_d:=[0, 1]^d$, that is
$$
  \Rc _d :=\left\{ \prod _{i=1}^d I_i  \,\,\, |\,\,\,  I_i =[a_i, b_i) \subset [0, 1] \r\}.
$$
%
%
The dispersion of a finite set $P\subset Q_d$ is defined as
$$
   \d (P) = \sup \{|B| \, \, \, | \, \, \, B\in \Rc _d, \, B\cap P=\emptyset\} .
$$
Then the minimal dispersion is defined as the function of two variables
--- the cardinality of a set of points $P\subset Q_d$ and the dimension, namely
$$
 \d^* (n, d) = \inf_{|P|=n} \d(P) .
$$
We also define its inverse as
$$
  N(\eps, d) = \min\{ n\in \N\,  | \,\,   \d^* (n, d)\leq \eps\} .
$$
Since in our proofs we use a random choice of points, it will be natural to
prove results in terms of the function $N(\eps, d)$
and then to provide the  corresponding (equivalent)  bounds for the
minimal dispersion itself.

\subsection {Known results.}

First we discuss the known bounds. In \cite{AHR} it was shown that for $\eps<1/4$,
\begin{equation}\label{first}
 (1-4\eps)\,  \frac{ \log_2 d}{4\eps } \leq  N(\eps, d) \leq \frac{2^{7d+1}}{\eps} ,
\end{equation}
where the upper bound is due to Larcher, improving the  previous bound via primorials
due to Rote and Tichy \cite{RT} (see also \cite{DJ})  and the lower bound is the first non-trivial
bound showing that the minimal dispersion grows with the dimension. Note that one trivially
has $\d^*(n, d) \geq 1/(n+1)$, hence  $N(\eps, d)\geq 1/\eps-1$.

Although estimates in (\ref{first}) are tight when the dimension $d$ is small and $\eps$ goes to 0,
 there is a huge gap between the upper and lower bounds when the dimension starts to grow
 to infinity. Using random choice of points uniformly distributed in $Q_d$, Rudolf \cite{Rud}
 obtained
\begin{equation}\label{rudbou}
   N(\eps, d) \leq \frac{8 d}{\eps}\log_2 \left(\frac{33}{\eps}\r)
\end{equation}
 (this bound with different numerical constants also follows from  much more general
 results in \cite{BEHW}, where the VC dimension
 of $\Rc_d$  was used, and from the fact that this VC dimension
 equals to $2d$). Estimate (\ref{rudbou}) is better than the upper bound in (\ref{first}) in the regime
$$
 \eps \geq \exp(-C^d),
$$
where $C>1$ is an absolute constant
(in this note we do not try to compute actual numerical values of absolute constants,
that is, constants independent of any other parameters, one can find them following the proofs).
Thus, if $\eps$ is not extremely small with respect to the dimension, the gap in bounds
is polynomial in $d$ and logarithmical in $1/\eps$. Another important feature of the Rudolf proof
is that a random choice of points uniformly distributed on $Q_d$ gives the result.

It was natural to conjecture that $N(\eps, d)$ behaves as $d/\eps$, especially in view
of corresponding bounds in the periodic setting (see below), however, surprisingly,
 Sosnovec \cite{Sos}
was able to improve the upper bound for $\eps <1/4$ to
$$
    N(\eps, d) \leq C_\eps \log_2 d ,
$$
where the order of magnitude of $C_\eps$ was essentially  $(1/\eps)^{(1/\eps)^2}$. This
dependence was significantly improved in \cite{UV} by Ullrich and Vyb\'iral, who showed that
$$
 C_\eps = \frac{2^7 }{\eps^2}\, \left(\log_2 \left(\frac{1}{\eps}\r)\r)^2
$$
works. They also conjectured that $N(\eps, d)$ behaves as $\log d/\eps$.
 The Sosnovec--Ullrich--Vyb\'iral upper bound is better in the regime
$$
  \eps \geq  \frac{C\, (\log_2 d)^2}{d}.
$$
The Sosnovec--Ullrich--Vyb\'iral proof is also based on a random choice of points,
but instead of the uniform distribution on $Q_d$ they use uniform distribution on a certain
lattice, gaining in the case of large $\eps$. We discuss this in more details  below.
Let us also mention that in the same paper Sosnovec proved that the function $N(\eps, d)$
completely changes the behaviour at $\eps=1/4$, more precisely, he proved that for every $\eps > 1/4$,
$$
  N(\eps, d)\leq 1 + \left\lf \frac{1}{\eps-1/4}  \r\rf.
$$
Thus, for $\eps > 1/4$, the function $N(\eps, d)$ is not growing with $d$.
Note that clearly $N(1/2, d) =1$ (by taking the point $(1/2, 1/2, ..., 1/2)$).
One can summarize the previously known upper bounds for $\eps\leq 1/4$ in
$$
   N(\eps, d)\leq
 \begin{cases}
     \frac{C\, \ln d }{\eps^2}\, \ln^2 \left(\frac{1}{\eps}\r), &\mbox{ if }\, \eps \geq  \frac{\ln^2 d}{d},\\
     \frac{C \, d }{\eps}\, \ln \left(\frac{1}{\eps}\r), &\mbox{ if }\, \frac{\ln^2 d}{d}\geq \eps \geq \exp(-C^d) ,\\
     \frac{C^d }{\eps}, &\mbox{ if }\, \eps \leq \exp(-C^d),
 \end{cases}
$$
where $1<C<1000$ is an absolute constant.

\subsection{New results}

In this note we improve the known bounds in the regime $\eps \geq \exp(-C^d)$.
Our first  result improves bounds when $\eps$ is not large.

\begin{thm}\label{th1}
There exists an absolute constant $C\geq 1$ such that the following holds.
Let $d\geq 2$ and $\eps \in (0, 1/2]$. Then
$$
 (i) \quad\quad\quad\quad     N(\eps, d)\leq \frac{C\, \ln d}{\eps }\,\, \, \, \ln \left(\frac{1}{\eps}\r),
 \quad\quad \mbox{ provided that } \quad  \eps \leq  \exp(-d),
 \quad\quad  \quad\quad  \quad\quad  \quad\quad
$$
$$
 (ii) \quad\quad \quad\,\,\,    N(\eps, d)\leq \frac{C \, d}{\eps } \,\, \ln  \ln \left(\frac{2}{\eps}\r),
 \quad\quad \,\,\, \mbox{ provided that } \quad  \eps \geq \exp(-d) .
 \quad\quad  \quad\quad  \quad\quad  \quad\quad
$$
Moreover,  the random  choice of points
with respect to the uniform distribution on the cube $Q_d$ gives
the result with high probability.
\end{thm}

We would like to emphasize  that if  $\eps \leq  \exp(-d)$ then, in view of (\ref{first}),
Theorem~\ref{th1} yields
$$
  \frac{\ln d}{6 \eps } \leq N(\eps, d)\leq \frac{C\, \ln d}{\eps }\, \, \ln \left(\frac{1}{\eps}\r),
$$
thus the gap in bounds is only logarithmical in $1/\eps$. In the second case
the improvement is only in substitution
of  $\ln(1/\eps)$ with $\ln \ln(1/\eps)$ comparing to Rudolf's bound.

Our proof is also based on a random choice of points. A standard way to use randomness
is to show that a certain ``good" event $\Ec$ holds with a non-zero probability. Equivalently,
one needs to show that the complement of $\Ec$, the event $\Ec^c$, holds with small probability.  In order to
do that, one tries to cover $\Ec^c$ by certain events, called  individual events,
  to obtain  good
bounds on probabilities of individual events, and then to use the union bound.
 In this scheme  one needs to have a good balance
between  (small) probabilities of  individual events and the (large but not too
large) size of the covering set. Since we need to prove that there exists a set $P$
of $n$ points such that there is no rectangle of volume $\eps$ without a point from $P$,
the natural idea would be to construct a finite set $\Nc$ of rectangles having reasonably large volume
and such that property {\it
``each rectangle in $\Nc$ contains a point from $P$"} implies the property {\it
``each rectangle in $\Rc_d$ of volume at least $\eps$ contains a point from $P$."}\,
In the case of uniform distribution on the cube $Q_d$, that is, in the case when the set  $P$ consists
of $N$ points independently drawn from the uniform distribution, an individual bound, that is, a bound
on the event that a given box $B\in \Nc$ contains a point from $P$, is simply given by the volume of $B$,
therefore the main difficulty is to construct the set $\Nc$ of not too large
cardinality. Rudolf used the concept of $\delta$-cover \cite{Rud, MG} to construct  $\Nc$ and to estimate
its size. We introduce the notion of $\delta$-net (see Definition~\ref{delta-net}), which fits better
for random procedure described above and allows to obtain better bounds on its size, see Propositions~\ref{net-gen}
and \ref{net-d}.


  As usual in probabilistic proofs, we obtain the result with high probability.
  Very recently, Hinrichs, Krieg, Kunsch, and Rudolf \cite{HKKR} investigated
the best bound that one can get using a random choice of points and showed that one cannot
expect anything better than
\begin{equation}\label{lrb}
  \max \left\{ \frac{c}{\eps} \ln\left(\frac{1}{\eps}\r),\, \frac{d}{2 \eps}\r\},
\end{equation}
where $c>0$ is an absolute constant. This in particular
shows that our bounds are almost best possible for this method (up to $\ln d$ in the first estimate and
up to $\ln \ln(1/\eps)$ in the second estimate).

In the case of large $\eps$ we can improve the bound. The next theorem provides better bounds
in the regime $\eps \geq (\ln^2 d)/(d\ln \ln (2d))$.

\begin{thm}\label{th2}
There exists an absolute constant $C\geq 1$ such that the following holds.
Let $d\geq 2$ and $\eps \in (0, 1/2]$ be such that $\eps \geq  \frac{\ln d}{d}$.  Then
$$
 N(\eps, d)\leq \frac{C \,\ln d}{\eps^2 }\, \, \ln \left(\frac{1}{\eps}\r).
$$
\end{thm}

This improves the Ullrich--Vyb\'iral bound by removing one $\ln(1/\eps)$ factor.
The proof of this theorem also
uses random points uniformly distributed on the cube $Q_d$, however, as
Hinrichs--Krieg--Kunsch--Rudolf's result shows, one cannot expect a bound better
than $d/\eps$, therefore one needs to adjust the distribution of the points.
One way to adjust randomness was suggested by  Sosnovec and then improved by
Ullrich and Vyb\'iral. They substituted the uniform distribution on the cube by
a uniform distribution on a certain lattice inside the cube. This led to the
logarithmic in $d$ upper bound (by the price of an additional factor $1/\eps$).
Careful analysis of their proofs in comparison with Rudolf's proof shows
that the main advantage of the use of a lattice is that the points on the lattice
are $\eps$-separated from the boundary of the cube. This leads to our adjustment
of the uniform distribution on the cube --- if a uniformly distributed over the cube
random point falls too close to the boundary
we slightly shift it to the interior, to ensure that it is $\eps$-separated from the boundary.
In the next section we introduce the function $\phi_\eps$, which serves this purpose.
Unfortunately, the size of $\delta$-nets is still too large, to deal with large $\eps$, so
we additionally introduce the notion of {\it dinets} --- {\it nets in the sense of dispersion}
(see Definition~\ref{defdinet}), which allows us to reduce the cardinality of a covering
set (see Proposition~\ref{dinet}) and hence to apply the union bound.

   The upper bonds for $\eps\leq 1/4$ from Theorems~\ref{th1} and \ref{th2} are summarized  in
$$
   N(\eps, d)\leq
 \begin{cases}
     \frac{C\, \ln d }{\eps^2}\, \ln \left(\frac{1}{\eps}\r), &\mbox{ if }\, \eps \geq  \frac{\ln^2 d}{d\ln \ln (2d)},\\
     \frac{C \, d }{\eps}\, \ln \ln \left(\frac{1}{\eps}\r), &\mbox{ if }\, \frac{\ln^2 d}{d\ln \ln (2 d)}\geq \eps \geq e^{-d} ,\\
     \frac{C \ln d }{\eps}\, \ln \left(\frac{1}{\eps}\r),  &\mbox{ if }\, e^{-d}\geq \eps \geq \exp(-C^d),\\
     \frac{C^d }{\eps}, &\mbox{ if }\, \eps \leq \exp(-C^d)
 \end{cases}
$$
or in the following picture showing the corresponding regions.

\bigskip

\bigskip

\includegraphics[width=0.8\textwidth]{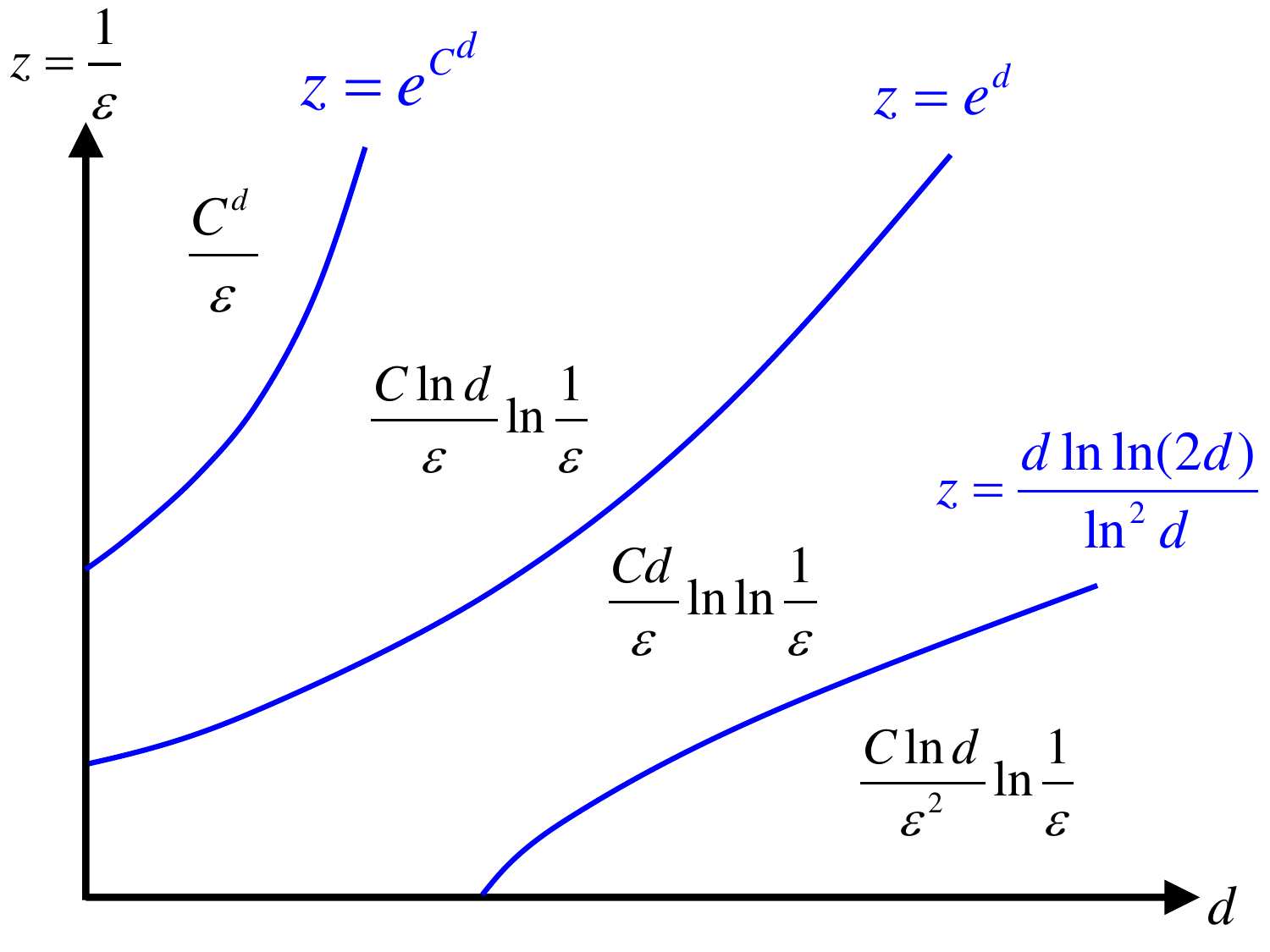}

\bigskip

In terms of the minimal dispersion, Theorems~\ref{th1} and \ref{th2} are equivalent to the following theorem.


\begin{thm}\label{disp1}
There exists an absolute constant $C\geq 1$ such that the following holds.
Let $d\geq 2$ and $n\geq 2\ln d$. Then
$$
 (i) \quad\quad     \d^*(n, d)\leq \frac{C\, \ln d}{n }\,\, \, \, \ln \left(\frac{n}{\ln d}\r),
 \quad\quad \mbox{ provided that } \quad\quad  n \geq   e^d d \ln d ,
 \quad\quad  \quad\quad  \quad\quad  \quad\quad
$$
$$
 (ii) \quad\,\,\,   \d^*(n, d)\leq  \frac{C \, d}{n} \,\, \ln  \ln \left(\frac{n}{d}\r),
 \quad\quad \,\, \mbox{ provided that } \quad \frac{d^2 \ln^ 2\ln d}{\ln^2 d}\leq n \leq e^d d \ln d,
 \quad\quad  \quad\quad
$$
 $$
 (iii) \quad\,  \d^*(n, d)\leq \left(\frac{C \,\ln d}{n }\, \, \ln \left(\frac{n}{\ln d}\r)\r)^{1/2},
 \quad\quad \mbox{ provided that } \quad\quad  n \leq  \frac{d^2 \ln^ 2\ln d}{\ln^2 d}.
 \quad\quad  \quad\quad  \quad\quad  \quad\quad
$$
Moreover, in the first two cases the random  choice of points
with respect to the uniform distribution on the cube $Q_d$ gives
the result with high probability.
\end{thm}

\subsection{Dispersion on the torus}

The corresponding  dispersion on the torus can be described in terms of
periodic axis parallel boxes. We denote such a set by $\wt \Rc _d$, that is
$$
  \wt \Rc _d :=\left\{ \prod _{i=1}^d I_i(a, b) \,\,\, |\,\,\,  a,  b\in  Q_d \r\},
$$
 where
$$
   I_i(a, b) :=
   \begin{cases}(a_i, b_i), &\mbox{ whenever }\, 0\leq a_i< b_i\leq 1 ,\\
    [0, 1]\setminus  [b_i, a_i], &\mbox{ whenever }\, 0\leq b_i< a_i\leq 1. \end{cases}
$$
The dispersion of a finite set $P\subset Q_d$ on the torus, the minimal dispersion on the torus,
and its inverse are defined in the same way as above, but using sets from $\wt \Rc _d$,
that is
$$
   \wt \d (P) = \sup \{|B| \, \, \, | \, \, \, B\in \wt \Rc _d, \, B\cap P=\emptyset\},
 \quad \quad  \quad
   \wt \d^* (n, d) = \sup_{|P|=n} \wt \d(P)  ,
$$
and
$$
 \wt N(\eps, d) = \min\{ n\in \N\,  | \,\,  \wt \d^* (n, d)\leq \eps\}  .
$$
 It is known that
$$
   \frac{d}{\eps} \leq  \wt N(\eps, d)\leq \frac{8 d}{\eps} \, \left(\ln d +  \ln \left(\frac{8}{\eps} \r)\r),
$$
where the lower bound was proved by Ullrich \cite{MU} and the upper bound is
due to Rudolf \cite{Rud} (since there are no good bounds on the VC dimension of
$\wt \Rc _d$, results of \cite{BEHW} are not directly applicable here).
We would like to emphasize that contrary to the non-periodic
case, even in the case of large $\eps$, the lower bound is at least $d$.
We improve the Rudolf upper bound in the case $\eps \leq 1/d$.

\begin{thm}\label{th-per}
There exists an absolute constant $C\geq 1$ such that the following holds.
Let $d\geq 2$ and $\eps \in (0, 1/2]$. Then
$$
 (i) \quad\quad\quad\quad    \wt N(\eps, d)\leq \frac{C\, \ln d}{\eps }\,\, \, \, \ln \left(\frac{1}{\eps}\r),
 \quad\quad \mbox{ provided that } \quad  \eps \leq  \exp(-d),
 \quad\quad  \quad\quad  \quad\quad  \quad\quad
$$
$$
 (ii) \quad\quad \quad\,\,\,  \wt  N(\eps, d)\leq \frac{C \, d\, \ln d}{\eps },  \quad\quad \quad
 \quad\quad \,\,\, \mbox{ provided that } \quad  \eps \geq \exp(-d).
 \quad\quad  \quad\quad  \quad\quad  \quad\quad
$$
Moreover,  the random  choice of points
with respect to the uniform distribution on the cube $Q_d$ gives
the result with high probability. Equivalently, for $d\geq 2$ and $n\geq 2d\ln d$ we have
$$
 (i) \quad\quad   \wt  \d^*(n, d)\leq \frac{C\, \ln d}{n }\,\, \, \, \ln \left(\frac{n}{\ln d}\r),
 \quad \mbox{ provided that } \quad\quad  n \geq   e^d d \ln d ,
 \quad\quad  \quad\quad  \quad\quad  \quad\quad
$$
$$
 (ii) \quad\,\,\, \wt  \d^*(n, d)\leq  \frac{C \, d\, \ln d}{n} ,
 \quad\quad \quad\quad\quad\,\, \mbox{ provided that }\quad \quad 2d\ln d \leq n \leq e^d d \ln d.
 \quad\quad  \quad\quad
$$
\end{thm}

Our bound on $\wt  N(\eps, d)$ reduces the factor $d$ in Rudolf's estimate to $\ln d$ in the case when
$\eps \leq \exp(-d)$ and removes the summand $\ln (1/\eps)$ if $ \exp(-d)< \eps \leq 1/d$. However,
if $\eps \geq 1/d$, it gives the same  order $(d \ln d)/\eps$.

The proof is the same as for Theorem~\ref{th1}, using random points and a $\delta$-net
constructed for periodic boxes. Unfortunately, in the construction of nets for the second
bound in Theorem~\ref{th1} and for the bound in Theorem~\ref{th2}, we essentially use
that boxes are not periodic and therefore the construction cannot be extended to the periodic case
(for  Theorem~\ref{th2} it is also clear in view of the Ullrich
lower bound on $\wt N(\eps, d)$). We would also like to note that the
Hinrichs--Krieg--Kunsch--Rudolf's result on best possible lower bound (\ref{lrb}) which may be obtained
by using random points uniformly distributed on the cube holds for the periodic setting as well,
therefore the factor $\ln (1/\eps)$ in our first estimate is unavoidable by this method.
In the second case, $\eps \geq \exp(-d)$, we have $\ln (1/\eps)\leq d$, so there is a hope
to remove $\ln d$ factor and to obtain the best possible estimate, on the other hand it is
possible that the bound is the best possible for this method.



\section{Nets, dinets, and a probability lemma}
\label{intro}

We need more notations. Given a positive integer $m$ we denote $[m]=\{1, 2, ..., m\}$.
Given $\eps >0$, we consider sets of (periodic) axis parallel of volume at least $\eps$,
$$
   \Be :=\Big\{ B\in \Rc _d   \,\,\, |\,\,\,  |B|\geq \eps \Big\}
   \qand  \Bep :=\Big\{ B\in \wt \Rc _d   \,\,\, |\,\,\,  |B|\geq \eps \Big\}.
$$

We introduce the following definition.

\begin{defin}[$\delta$-net for $\Be$] \label{delta-net}
 Given $\eps, \delta >0$ we say that $\Nc \subset \Rc _d$
is a $\delta$-net for $\Be$ if for every $B\in \Be$
there exists $B_0\in \Nc$ such that $B_0\subset B$ and
$$
 |B_0|\geq (1-\delta) |B|.
$$
We define a $\delta$-net for $\Bep$ in a similar way.
\end{defin}

To deal with large $\eps$ with respect to the dimension, say when $\eps\geq 1/d$,
we adjust the definition of a $\delta$-net by   introducing the notion of {\it $\delta$-dinet} ---
{\it a $\delta$-net in a sense of dispersion}.
The key idea leading to this approach is an observation that we do not need to consider
points which are too close to the boundary of the cube $Q_d$. As we mentioned in the introduction,
this idea was already implicitly used in \cite{Sos, UV}.
First given $\eps\in (0, 1/2)$ define an auxiliary function
$\phi_\eps \, :\, [0,1]\to [\eps, 1-\eps]$ by
$$
    \phi_\eps (t) = \begin{cases} \eps &\mbox{ if  }\, 0\leq t< \eps ,\\
       t  &\mbox{ if  }\, \eps\leq t\leq 1-\eps ,\\
       1-\eps  &\mbox{ if  }\, 1-\eps< t\leq 1
    . \end{cases}
$$
Given $x\in Q_d$ we also write $\phi_\eps(x)$ for $\{\phi_\eps (x_i)\}_{i=1}^d$.

\begin{defin}[$\delta$-dinet for $\Be$] \label{defdinet}
 Given $\eps, \delta >0$ we say that $\Nc \subset \Rc _d$
is a $\delta$-dinet for $\Be$ if for every $B\in \Be$
there exists $B_0\in \Nc$ such that
$$
 |B_0|\geq (1-\delta) |B|
$$
and such that for every $x\in  Q_d$ the following implication holds
$$
   x\in B_0 \quad \Longrightarrow \quad  \phi_\eps(x) \in B
$$
\end{defin}

Note that every $\delta$-dinet $\Nc$ for $\Be$ has the following property
allowing to bound from above the number of points needed to
have a given dispersion, namely, for every $n\geq 1$, every set
of points $P=\{x_1, ..., x_n\}\subset Q_d$ the statement
{\it ``each box from $\Nc$ contains at least one point from $P$''}
implies the statement
{\it ``each box from $\Be$ contains at least one point from $\phi_\eps(P)$''}.}

A variant of the following lemma using random points and the union bound
was proved in  \cite{Rud} (see Theorem~1 there). We provide a proof for
completeness.

\begin{lemma}\label{unb}
Let $d\geq 1$ and  $\eps, \delta \in (0,1)$. Let $\Nc$ be either a $\delta$-net for $\Be$
or a $\delta$-dinet for $\Be$ and
let   $\wt \Nc$ be a $\delta$-net for $\Bep$. Assume both $|\Nc|\geq 3$ and $|\wt \Nc|\geq 3$.
 Then
$$
   N(\eps, d) \leq  \frac{3\ln |\Nc|}{(1-\delta) \eps}
    \qand  \wt N(\eps, d)\leq \frac{3\ln |\wt \Nc|}{(1-\delta) \eps} .
$$
\end{lemma}

\begin{rem}\label{reml}
As usual for proofs involving the union bound,
 our proof shows that the random choice of $N=\lfloor \frac{3\ln |\Nc|}{(1-\delta) \eps}\rfloor$ 
 (resp. $N=\lfloor \frac{3\ln |\wt \Nc|}{(1-\delta) \eps}\rfloor$) points gives the result with high
probability, more precisely with probability at least $1-1/|\Nc|$.
In the case of  $\delta$-nets the randomness is with respect to the independent
uniform choice of points on $Q_d$, while in the case of $\delta$-dinets one needs
to adjust the choice of independent uniformly distributed points by the function
$\phi_\eps$.
\end{rem}

\begin{proof}
We show a proof for a $\delta$-net $\Nc$ for $\Be$, the other two cases are  the same.
Let $\Nc$ be a $\delta$-net for $\Be$. Consider $N$ independent random points
$X_1$, ..., $X_N$ uniformly chosen from $Q_d$. By the definition of a $\delta$-net, it is enough
to show that for every $B\in \Nc$ with $|B|\geq v:=(1-\delta) \eps$ there exists $j\leq N$
such that $X_j\in B$. Fix such a box $B$. Using that the volume of $B$ is at least $v$
and the independence of $X_j$'s, we obtain
$$
  \Prob \left( \left\{\forall j\leq N :\, \, \, X_j\notin B \r\}\r) \leq  \left(1 - v\r)^N < \exp(-vN).
$$
Therefore, by the union bound,
$$
  \Prob \left( \left\{\exists B\in \Nc :   \, \, \,|B|\geq v \qand \forall j\leq N :\, \, \, X_j\notin B \r\}\r)
  <
   |\Nc| \exp(-vN).
$$
Thus, as far as $|\Nc| \exp(-vN)\leq 1$, there exists a realization of $X_j$'s with the desired
property. Moreover, if $N=\lfloor \frac{3\ln |\Nc|}{(1-\delta) \eps}\rfloor$ then the ``good" 
probability is $1-|\Nc| \exp(-vN)\geq 1 - 1/ |\Nc| $. 
This implies both Lemma~\ref{unb} and Remark~\ref{reml}.
\end{proof}

\section{Cardinality of nets}
\label{card}

As is seen from Lemma~\ref{unb} and Remark~\ref{reml}, to prove  our theorems,
it is enough to construct nets of not so large cardinality.
The next simple observation is one of key ideas in our estimates.
Let $\eps >0$ and let $\ell_1, ..., \ell _d>0$ be such that
$$
   \prod _{i=1}^{d} \ell_i \geq \eps.
$$
Denote by $\sigma = \sigma (\ell_1, ..., \ell _d)$ a permutation
such that
\begin{equation}\label{perm}
     \ell _{\sigma(1)} \leq  \ell _{\sigma(2)} \leq ...  \leq  \ell _{\sigma(d)}
\end{equation}
(for each sequence we fix one such permutation).
Then for every $j\geq 1$ we clearly have
\begin{equation}\label{lbound}
   \ell _{\sigma(j)} \geq \left(\prod _{i=1}^{j} \ell_{\sigma (i)}\r)^{1/j} \geq \eps ^{1/j}
   > 1-\frac{\ln(1/\eps)}{j}.
\end{equation}
 A naive approach to approximate rectangles from $\Be$ is to say that given a
 rectangle $B=\prod _{i=1}^{d} I_i \in \Be$  the smallest
 length $\ell_i=|I_i|$ is at least $\eps$. Therefore, we can take
 $(1/(4\eps))$-net $\Mc$ in $[0,1]$ and approximate each $I_i$ with segments
 having endpoints in $\Mc$. This approach would lead to a net of the order
 $(1/(4\eps))^{2d}$, which is not acceptable for
 our purpose (this would also lead to a huge loss in volume, but
 already the size of a net is too large). Instead, we use formula (\ref{lbound}),
 to say that the larger $i$ the coarser
  net in $[0,1]$ is needed in order to approximate the corresponding
 interval $I_{\sigma(i)}$. Of course, simultaneously, we need to control the loss
 in volume in our approximation. The next proposition utilizes this idea.
 It works for both the periodic and  non-periodic settings.
 Since we will be using this result in several dimensions,   it would be convenient
 to formulate it for boxes in $\R^m$.

\begin{prop}\label{net-gen}
Let $m\geq 2$ be an integer and $\eps \in (0,1)$.
There are $(1/2)$-nets $\Nc$ and $\wt \Nc$  for $\Bem$ and $\Bepm$
respectively, each of them of cardinality at most
$$
    \frac{(14 m)^{4m}}{\eps ^{2\log _2(2 m)}}.
$$
\end{prop}

\begin{rem}\label{net-rem}
If $m=2^k$ for some integer $k$ then our proof gives slightly better estimate, namely
$$\frac{(24 m)^{2m}}{\eps ^{2\log _2 m}}.$$
\end{rem}

\begin{rem}\label{proof1}
 Clearly, Lemma~\ref{unb} and Remark~\ref{reml} combined with this proposition
 (applied with $m=d$) yield
 Theorem~\ref{th-per} as well as the first bound in Theorem~\ref{th1}.
\end{rem}

\begin{proof}
The construction of nets in $\Bem$ and $\Bepm$ are essentially the same.
 We provide a proof for a net in $\Bepm$, since  the proof for a net in $\Bem$
is somewhat easier --- we do not need to consider intervals $I_i(a, b)$ with $a_i>b_i$.

 Fix $k\geq 1$ such that $2^k\leq m<2^{k+1}$. Fix a partition
 of $[m]$  into $k+1$ disjoint sets $A_1$, ..., $A_{k+1}$ with
 $|A_1|=2$, $|A_{k+1}|=m-2^k$ (this set is empty if $m=2^k$),
 and $|A_j|=2^{j-1}$ for $2\leq j\leq k$.
 For $j\leq k+1$ denote
 $$
 \delta^{(j)} =2^{-k-3} \eps^{2^{1-j}} \qand
  D_{j}=\{0, \delta^{(j)}, 2\delta^{(j)}, ..., s_j \delta^{(j)}\},
 $$
 where $s_j=\lf 1/\delta^{(j)}\rf$ (note that dealing with $\Bem$ we do not need to
have $0$ in $D_{j}$).

 We are now ready to define a part of our net
 corresponding to this partition of $[m]$ as the set
%
$$
 \Nc _*(A_1, ..., A_{k+1}) :=
  \Big\{  \prod _{i=1}^m I_i(x, y) \,\,\, | \,\,\, x, y \in Q_m,\,\, \forall j\leq k+1 \,\, \forall \, i\in A_j: \, \, x_i\ne y_i\in  D_{j}
  \Big\}.
$$
Then the cardinality of this set  can be estimated as
\begin{align*}
  |\Nc _*(A_1, ..., A_{k+1})|
  &\leq  \prod _{j=1}^{k+1} \prod _{i\in A_j}   |D_{j}|( |D_{j}|-1)  \leq
    \prod _{j=1}^{k+1} \prod _{i\in A_j}\frac{ 2}{ \left(\delta^{(j)}\r)^2} \leq 2^m
    \prod _{j=1}^{k+1} \prod _{i\in A_j} \frac{ 4^{k+3}}{\eps^{2^{2-j}}}
   \\ &\leq 2^m
      \left(\frac{ 4^{k+3}}{\eps^{2}} \r)^{2} \,\,
      \prod _{j=2}^{k+1}  \left(\frac{ 4^{k+3}}{\eps^{2^{2-j}}} \r)^{2^{j-1}}= 2^m
      \frac{ 4^{2^{k+1}(k+3)} }{ \eps^{2(k+1)} }\leq \frac{2^m (64 m^2)^{2m}}{\eps ^{2\log _2(2 m)}}
\end{align*}
(note that if $m=2^k$, then the set $A_{k+1}$ is empty and $j$ runs between 1 and $k$,
which leads to the bound from Remark~\ref{net-rem}).

To complete the construction, we take the union over all partitions of $[m]$ into such sets $A_1, ..., A_{k+1}$,
$$
   \Nc :=\bigcup \Nc _*(A_1, ..., A_{k+1}) .
$$
The number of partitions can be estimated as
$$
  {m\choose 2^k} {2^k\choose 2^{k-1}} {2^{k-1}\choose 2^{k-2}} ... {4\choose 2} \leq 2^{2m},
$$
hence
$$
  |\Nc|\leq \frac{2^{3m} (8 m)^{4m}}{\eps ^{2\log _2(2 m)}}\leq  \frac{ (14 m)^{4m}}{\eps ^{2\log _2(2 m)}}.
$$

 It remains to show that $\Nc$ is indeed a $(1/2)$-net for $\Bepm$. Let $a, b\in Q_m$ and
 $B=\prod _{i=1}^{m} I_i(a, b)$ be of volume at least $\eps$. For $i\leq m$
 let $\ell_i$ be the length of  $I_i(a, b)$. Let $\sigma = \sigma (\ell_1, ..., \ell _m)$
  be the  permutation defined by (\ref{perm}).
 Consider the following partitions of $[m]$,
 $$
   A_1^\sigma=\sigma(\{1, 2\}), \,\,\, A_{k+1}^\sigma =\sigma(\{2^k + 1, ..., m\}), \qand
   A_j^\sigma=\sigma(\{2^{j-1} +1, ..., 2^j\}),
 $$
 $2\leq j\leq k$,
 and note that by (\ref{lbound}) for every $j\leq k+1$ and every $i\in A_j^\sigma$ one has
 \begin{equation}\label{len}
    |I_i(a, b)|=\ell_i \geq \ell_{\sigma(2^{j-1})}\geq \eps ^{1/2^{j-1}} = 2^{k+3}\delta^{(j)} .
 \end{equation}
Take a box $B_0=\prod _{i=1}^{m} I_i(x, y)$ from $\Nc _*(A_1^\sigma, ..., A_{k+1}^\sigma)$
such that for every $j\leq k+1$ and every $i\in A_j^\sigma$ one has
$$
    a_i \leq x_i , \,\,\,\, b_i\geq y_i,\,\,\,\, x_i-a_i \leq
    \delta ^{(j)}, \qand  b_i-y_i \leq \delta ^{(j)}
$$
 (if $a_i > s_j \delta^{(j)}$ we take $x_i=0$).
 The lower bound (\ref{len}) on the length of $I_i(a, b)$ implies that
 $I_i(x, y)\subset I_i(a, b)$. Thus, $B_0\subset B$ and, using (\ref{len}) again,
 \begin{align*}
   |B_0| &= \prod _{j=1}^{k+1} \prod _{i\in A_j^\sigma}   |I_i(x, y)| \geq
    \prod _{j=1}^{k+1} \prod _{i\in A_j^\sigma} \left(\ell_i-2\delta^{(j)}  \r)=
    \prod_{i=1}^{m} \ell_i \,\,
    \prod _{j=1}^{k+1} \prod _{i\in A_j^\sigma} \left(1-\frac{2\delta^{(j)}}{\ell_i}  \r)
   \\ &\geq
   |B| \, \prod _{j=1}^{k+1}  \left(1-\frac{2\delta^{(j)}}{\ell_{\sigma(2^{j-1})}}  \r)^{|A_j^\sigma|}
    \geq |B|  \prod _{j=1}^{k+1}  \left(1-\frac{1}{2^{k+2}}  \r)^{|A_j^\sigma|}\geq
    |B|  \left(1-\frac{1}{2^{k+2}}  \r)^{m}
    \\ &\geq
    |B|  \left(1-\frac{1}{2 m}  \r)^{m}\geq \frac{1}{2}\, |B| .
\end{align*}
 This completes the proof.
\end{proof}

Next we show how to improve the bound of Proposition~\ref{net-gen} for non-periodic
boxes in the case when $\eps$ is not very small with respect to dimension, say, when
$4\ln (1/\eps) \leq d$. The key observation here is that in the case $4\ln (1/\eps) \leq d$
a rectangle  $B=\prod _{i=1}^{d} I_i \in \Be$ has many intervals $I_i$ of length
close to one, namely, by (\ref{lbound}), $|I_{\sigma(i)}| \geq 1- 1/L$ whenever
$i\geq L \ln (1/\eps)$. For such an interval we do not need to take
a net in  $[0,1]$ in  order to approximate the end points --- it is enough
to approximate the left end point by a net in $[0, 1/L]$ and the right
 end point by a net in $[1-1/L, 1]$. This leads to a significant improvement in
 the size of the net. Of course, this approach cannot work for periodic boxes.

\begin{prop}\label{net-d}
Let $d\geq 4$ be an integer, $\eps \in (0, 1/4]$ and assume that $d\geq 4\ln (1/\eps)$.
Then $\Be$  admits a $(3/4)$-net
of cardinality at most
$$
    \exp\left(C d \ln \ln (1/\eps) \r),
$$
where $C\geq 1$ is an absolute constant.
\end{prop}

\begin{rem}\label{proof1}
 Clearly, Lemma~\ref{unb}  and Remark~\ref{reml} combined with this proposition  yield
  the second bound in Theorem~\ref{th1}.
\end{rem}

\begin{proof}
The proof is similar to the proof of Proposition~\ref{net-gen}, but
we deal more carefully with the approximation of long segments.

Set $k$ to be the smallest integer such that $2^k\geq 2\ln (1/\eps)$ and let $m=2^k$.
Clearly, $k\geq 1$, $m\geq 2$.  Then $d\geq 4\ln (1/\eps)>m$.
Fix an integer  $n\geq k$ such that $2^n\leq d<2^{n+1}$.
 Fix a partition of $[d]$  into $n-k+2$ disjoint sets $A_0$, ..., $A_{n-k+1}$ with
 $|A_0|=m$,  $|A_{n-k+1}|=d-2^n$ (this set is empty if $d=2^n$),
 and $|A_j|=2^{k+j-1}$ for $1\leq j\leq n-k$.
 Denote
 $$
 \delta =\frac{1}{8 d}, \quad
   D_1=\{\delta, 2\delta, ..., s \delta\}, \qand
    D_2=\{1-\delta, 1-2\delta, ..., 1- s \delta\},
 $$
where $s_j=\lf 1/\delta\rf$.

Next, for  every $1\leq j\leq n-k+1$ we consider the set $P_j\subset D_1\times D_2$ of all pairs $(p, q)$
satisfying $p\in D_1$, $q\in D_2$, $p<q$, and
$$
    p \leq 2^{1-k-j}\ln(1/\eps) +\delta \qand
     q\geq  1- 2^{1-k-j}\ln(1/\eps) -\delta.
$$
Using $2^k\geq 2\ln (1/\eps)$ and $\delta=1/(8d)$,
we observe that the cardinality of $P_j$ is
$$
   |P_j|\leq  \left(\frac{2^{1-k-j}\ln(1/\eps)}{\delta}+1\r)^2\leq
    \left(\frac{8d}{2^j}+1\r)^2\leq \frac{d^2}{4^{j-2}}
$$
Let $\Nc_0 (A_0)$ be the $(1/2)$-net of cardinality at most
$$
 n_0:=\frac{(24 m)^{2m}}{\eps ^{2\log _2 m}}
$$
for $\Bem$ from Proposition~\ref{net-gen}
constructed in $\R^{A_0}$
(see also Remark~\ref{net-rem}).
Let
$\Nc _*=\Nc _*(A_0, ..., A_{n-k+1})$ be the set of all boxes
$
    \prod _{i=1}^d [x_i, y_i)
$
such that
$$
 \prod _{i\in A_0} [x_i, y_i)\in \Nc_0 (A_0)
$$
 and for every
$1\leq j\leq n-k+1$ and for every $i\in A_j$ the pair $(x_i, y_i)\in P_j$.
Then, using $2^n\leq d<2^{n+1}$ and $m=2^k$,  the cardinality of $\Nc _*$ can be estimated as
\begin{align*}
  |\Nc _*|
  &\leq |\Nc_0 (A_0)|\,  \prod _{j=1}^{n-k+1} \prod _{i\in A_j}   |P_{j}| \leq
  n_0\, \prod _{j=1}^{n-k+1}
     \left(\frac{d^2}{4^{j-2}}\r)^{|A_j|} \leq
     n_0\, \prod _{j=1}^{n-k+1}
     \frac{d^{2^{k+j}}}{4^{(j-2) 2^{k+j-1}}}
   \\ &\leq
     n_0 \, \frac{ d^{2^{n+2}}}{4^{(n-k-3) 2^{n+1}}}
    = n_0 \,\left(\frac{ 4^{k+4}d^{2}}{4^{n+1}} \r)^{2^{n+1}}
    \leq  n_0 \,\left( 4^4 m^2 \r)^{2d} = \frac{(24 m)^{2m}}{\eps ^{2\log _2 m}}\,\,\left( 16 m \r)^{4d}.
\end{align*}

Finally we define our net as  the union over all partitions of $[d]$ into such sets
$A_0, ..., A_{n-k+1}$,
$$
   \Nc :=\bigcup \Nc _*(A_0, ..., A_{n-k+1}) .
$$
The number of partitions can be estimated as
$$
  {d\choose 2^n} {2^n\choose 2^{n-1}} {2^{n-1}\choose 2^{n-2}} ... {2^{k+1}\choose 2^k} \leq
  2^{2d-2^{k+1}}\leq  2^{2d-2m},
$$
hence
$$
  |\Nc|\leq  2^{2d-2m}\, \frac{(24 m)^{2m}}{\eps ^{2\log _2 m}}\,\,\left( 16 m \r)^{4d}
  \leq  \frac{(12 m)^{2m}}{\eps ^{2\log _2 m}}\, \left( 24 m \r)^{4d} \leq
  \frac{(24 m)^{6d}}{\eps ^{2\log _2 m}}.
$$
Using that $m\leq 4\ln(1/\eps)\leq d$, we obtain
$$
 |\Nc|\leq \exp\left(6d \ln (24m) + 2(\log_2 m)(\ln(1/\eps) \r)\leq  \exp\left(C d \, \ln \ln (1/\eps) \r),
$$
where $C\geq 1$ is an absolute constant.

 It remains to show that $\Nc$ is indeed a $(3/4)$-net for $\Be$. Let $a, b\in Q_d$ with $a_i<b_i$ for all $i\leq d$,
 and $B=\prod _{i=1}^{d} [a_i, b_i)$ be of volume at least $\eps$. For $i\leq d$
 let   $\ell_i=b_i-a_i$.
 Let $\sigma = \sigma (\ell_1, ..., \ell _d)$
  be the  permutation defined by (\ref{perm}).
 Consider the following partitions of $[d]$,
 $$
   A_0^\sigma=\sigma([m]), \,\,\, A_{n-k+1}^\sigma =\sigma(\{2^n + 1, ..., d\}), \qand
   A_j^\sigma=\sigma(\{2^{k+j-1} +1, ..., 2^{k+j}\}),
 $$
 $1\leq j\leq n-k$.

 Fix for a moment
 $1\leq j\leq n-k+1$ and $i\in A_j^\sigma$.
 Using $2^k\geq 2\ln(1/\eps)$ and (\ref{lbound}), we observe  that
 \begin{equation}\label{len-2}
    b_i-a_i=\ell_i \geq \ell_{\sigma(2^{k+j-1})}> 1-\frac{\ln(1/\eps)}{2^{k+j-1}}\geq 1-2^{-j}.
 \end{equation}
 Take a pair $(x_i, y_i) \in D_1\times D_2$ satisfying
 $$
    a_i \leq x_i , \,\,\,\, b_i\geq y_i,\,\,\,\, x_i-a_i \leq
    \delta, \qand  b_i-y_i \leq \delta .
$$
 Then $y_i-x_i\geq b_i-a_i-2\delta > 1-2^{-j}- 2\delta > 0$ and
$$
  y_i\geq b_i-\delta >  1-\frac{\ln(1/\eps)}{2^{k+j-1}} -\delta \qand
  x_i\leq a_i+\delta < \frac{\ln(1/\eps)}{2^{k+j-1}} + \delta ,
$$
in other words the pair $(x_i, y_i)\in P_j$.

Consider the box
 $B_0=\prod _{i=1}^{d} [x_i, y_i)$
such that for every $1\leq j\leq n-k+1$ and every $i\in A_j^\sigma$ the pair
$(x_i, y_i)$ is constructed as above and
where
$$
  B_{0}'=\prod _{i\in  A_0^\sigma} [x_i, y_i)\in \Nc_0 ( A_0^\sigma)\quad \quad \mbox{ approximates }
 \quad \quad B'=\prod _{i\in  A_0^\sigma} [a_i, b_i)
$$
as in Proposition~\ref{net-gen} (note that $m$-dimensional volume of $B'$ is at least
$\eps$, so $B'\in \Bem$).
Then by construction $B_0 \in \Nc _*(A_1^\sigma, ..., A_{k+1}^\sigma)$, $B_0\subset B$, and
$$
   | B_{0}'|\geq \frac{1}{2}\, |B'| = \frac{1}{2}\, \prod _{i\in  A_0^\sigma} \ell _i.
$$
Furthermore, using $\delta=1/(8d)$ and the bound (\ref{len-2})  again,
 \begin{align*}
   |B_0|&=|B_0'|\, \prod _{j=1}^{n-k+1} \prod _{i\in A_j^\sigma}   |y_i- x_i| \geq
    \frac{1}{2}\, \prod _{i\in  A_0^\sigma} \ell _i\, \prod _{j=1}^{n-k+1} \prod _{i\in A_j^\sigma} \left(\ell_i-2\delta  \r)
    \\&=    \frac{1}{2}\, \prod_{i=1}^{d} \ell_i \,\,
    \prod _{j=1}^{n-k+1} \prod _{i\in A_j^\sigma} \left(1-\frac{2\delta}{\ell_i}  \r)
  \geq 
   \frac{1}{2}\, |B|\,  \prod _{j=1}^{n-k+1}  \left(1-\frac{2\delta}{1-2^{-j}}  \r)^{|A_j^\sigma|}
   \\& \geq
    \frac{1}{2}\, |B|\,  \left(1-\frac{1}{2d}  \r)^{d-2^n}  \,\, \prod _{j=1}^{n-k}  \left(1-\frac{1}{2d}  \r)^{2^{k+j-1}}\geq
     \frac{1}{2}\, |B|\,  \left(1-\frac{1}{2 d}  \r)^{d}\geq  \frac{1}{4}\, |B|.
\end{align*}
 This completes the proof.
\end{proof}

Finally, we want to improve bounds in the case of large $\eps$.
The following proposition is an almost immediate consequence of Proposition~\ref{net-gen}
and definitions. The key observation here is also the fact that  a rectangle
$B=\prod _{i=1}^{d} I_i \in \Be$ has many intervals $I_i$ of the length
close to one, but now they will be so close to one, that we can substitute them
just by $[0, 1]$. More precisely, using our function $\phi_\eps$, if the length
of $I_i$ is at least $1-\eps$ then for every $z\in [0, 1]$ one has $\phi_\eps(z)\in I_i$,
hence we do not need to approximate such intervals. This  leads to our
definition of a dinet and to better bounds of cardinality of  dinets versus regular
nets. Unfortunately, this also leads to an additional factor $1/\eps$ in the final bound.
As in the previous proposition, such an approach essentially uses that we are in
the non-periodic setting.

\begin{prop}\label{dinet}
Let $d\geq 4$ be an integer, $\eps \in (0, 1/2]$ and assume that $d\geq (\ln (1/\eps))/\eps$.
There is a $(1/2)$-dinet $\Nc$ for $\Be$  of cardinality at most
$$
    \exp\left( \frac{9  \ln(1/\eps) \ln (18d) }{\eps}\r) .
$$
\end{prop}

\begin{rem}
 Clearly, Lemma~\ref{unb}  combined with this proposition  yields
   Theorem~\ref{th2}. We can also use  Remark~\ref{reml} to claim
   that a random choice of points works with high probability, but here the randomness will be
   with respect to the uniform distribution on the cube adjusted by the function
   $\phi_\eps.$
\end{rem}

\begin{proof}
Fix the smallest integer $m\geq (\ln(1/\eps) )/\eps$.
Given subset $A\subset[d]$ of cardinality $m$, let
 $\Nc_0 (A)$ be the $(1/2)$-net of cardinality at most
$$
 n_0:=\frac{(14 m)^{4m}}{\eps ^{2\log _2 (2m)}}
$$
for $\Bem$ from Proposition~\ref{net-gen}
constructed in $\R^{A}$. Let
$\Nc _*(A)$ be the set of all boxes
$
    \prod _{i=1}^d [x_i, y_i)
$
such that
$$
 \prod _{i\in A} [x_i, y_i)\in \Nc_0 (A)
$$
 and for every $i\notin A$, $[x_i, y_i)=[0,1)$.  Let
 $$
    \Nc = \bigcup_{A\subset [d]\atop |A|=m}  \Nc_*(A).
 $$
Then the cardinality of $\Nc$ is at most
$$
  {d\choose m } n_0 \leq \left(\frac{ed}{m}\r)^m n_0\leq
  \frac{(14^{4} e   m^3 d)^{m}}{\eps ^{2\log _2 (2m)}} \leq
  \frac{(18  d)^{4 m}}{\eps ^{2\log _2 (2m)}}.
$$
Since $(\ln(1/\eps) )/\eps\leq m \leq d$ and $m\leq 2(\ln(1/\eps) )/\eps$,
 this implies
$$
  |\Nc|\leq \exp( 4 m \ln (18d) + 2\log _2 (2m)  \ln(1/\eps))
  \leq \exp\left( \frac{9  \ln(1/\eps) \ln (18d) }{\eps}\r) .
$$

Now we show that $\Nc$ is a $(1/2)$-dinet  for $\Be$.
Let $a, b\in Q_d$ with $a_i<b_i$ for all $i\leq d$,
 and $B=\prod _{i=1}^{d} [a_i, b_i)$ be of volume at least $\eps$. For $i\leq d$
 let   $\ell_i=b_i-a_i$.
 Let $\sigma = \sigma (\ell_1, ..., \ell _d)$
  be the  permutation defined by (\ref{perm}) and
  denote $A^\sigma=\sigma([m])$.
Consider the box
 $B_0=\prod _{i=1}^{d} [x_i, y_i)$
such that $[x_i, y_i)=[0, 1)$ for every $i\notin A^\sigma$ and
$$
  B_{0}'=\prod _{i\in A^\sigma} [x_i, y_i)\in \Nc_0 ( A^\sigma)\quad \quad \mbox{ approximates }
 \quad \quad B'=\prod _{i\in  A^\sigma} [a_i, b_i)
$$
as in Proposition~\ref{net-gen} (note that $m$-dimensional volume of $B'$ is at least
$\eps$, so $B'\in \Bem$).
Then by construction $B_0 \in \Nc _*(A^\sigma)$, and
$$
  |B_0| =  | B_{0}'|\geq \frac{1}{2}\, |B'| = \frac{1}{2}\, \prod _{i\in  A^\sigma} \ell _i
  \geq \frac{1}{2}\, |B|.
$$

Finally assume that $z\in B_0$. If $i\notin  A^\sigma$ then using
(\ref{lbound}) and $m\geq (\ln(1/\eps) )/\eps$ we have
$$
 b_i-a_i = \ell_i \geq \ell_{\sigma (m)} >1-\eps.
$$
Therefore,
$\phi_\eps (z_i) \in [\eps, 1-\eps] \subset [a_i, b_i)$.
Assume $i\in A^\sigma$. Note that in this case
$$
 z_i \in [x_i, y_i)\subset [a_i, b_i),
$$
and $b_i-a_i\geq \eps$ (otherwise $|B|<\eps$).
If $\eps \leq z_i \leq 1-\eps$ then $\phi_\eps (z_i) = z_i$ hence $\phi_\eps (z_i)  \in  [a_i, b_i).$
If $0\leq z_i < \eps$ then the interval $[a_i, b_i)$ contains a point smaller than
$\eps$ and has length at least  $\eps$.  Then it must contain $\eps= \phi_\eps (z_i)$.
Similarly, if  $1-\eps< z_i < 1$ then $[a_i, b_i)$ must contain $1-\eps= \phi_\eps (z_i)$.
This proves that if $z\in B_0$ then $\phi_\eps (z) \in B$. Thus, $\Nc$ is a $(1/2)$-dinet
for $\Be$.
This completes the proof.
\end{proof}


\subsection*{Acknowledgments}
The author was introduced to this problem during the 2017 MFO workshop
{\it ``Perspectives in High-dimensional Probability and Convexity."}\,
The author is grateful to MFO,  to the organizers, and participants
of the workshop.
The author is also grateful to A.~Zelnikov for his help with the picture.

\address

\end{document}